\documentclass[a4paper, reqno]{amsart}

\usepackage{amsmath, amsthm, amssymb}
\usepackage{color}
\usepackage{fixltx2e}
\usepackage{enumitem}
\usepackage{lmodern}

\theoremstyle{plain}
\newtheorem{thm}{Theorem}[section]
\newtheorem{lemma}[thm]{Lemma}
\numberwithin{equation}{section}

\let\originalleft\left
\let\originalright\right
\renewcommand{\left}{\mathopen{}\mathclose\bgroup\originalleft}
\renewcommand{\right}{\aftergroup\egroup\originalright}
\renewcommand\Re{\operatorname{Re}}

\newcommand\twoln[2]{{\substack{#1 \\ #2}}}
\newcommand\thrln[3]{{\substack{#1 \\ #2 \\ #3}}}
\newcommand\Res[1]{\underset{#1}{\operatorname{Res}}}
\newcommand\smod[1]{\, (#1)}
\newcommand\BigO[1]{\mathcal{O}\left(#1\right)}
\newcommand\supp{\operatorname{supp}}

\begin{document}

\title{The shifted convolution of generalized divisor functions}

\author{Berke Topacogullari}
\address{Mathematisches Institut, Bunsenstra{\ss}e 3-5, D-37073 G\"ottingen, Germany}
\email{btopaco@uni-goettingen.de}
\thanks{This work was supported by the Volkswagen Foundation.}

\bibliographystyle{abbrv}

\begin{abstract}
  We prove an asymptotic formula for the shifted convolution of the divisor functions $d_k(n)$ and $d(n)$ with $k \geq 4$, which is uniform in the shift parameter and which has a power-saving error term, improving results obtained previously by Fouvry and Tenenbaum and, more recently, by Drappeau.
\end{abstract}

\maketitle

\section{Introduction}

The problem we are concerned with in this paper is the asymptotic evaluation of the following shifted convolution sums,
\[ D_k^+(x, h) := \sum_{n \leq x} d_k(n) d(n + h) \quad \text{and} \quad D_k^-(x, h) := \sum_{n \leq x} d_k(n + h) d(n), \]
where \(h\) is a positive integer, and where \( d_k(n) \) stands for the number of ways to write \(n\) as a product of \(k\) positive integers.

For \( k = 2 \), this is a classical problem also known as the binary additive divisor problem, and has been studied in the past decades by many authors.
As an example, it is known that
\begin{align} \label{eqn: result for D_2(x, h)}
  D_2^\pm(x, h) = x P_{2, h}(\log x) + \BigO{ x^{\frac23 + \varepsilon}} \quad \text{for} \quad h \ll x^\frac23,
\end{align}
with \( P_{2, h} \) a quadratic polynomial depending on \(h\), a result we have cited from Motohashi \cite{Mot94}, where a detailed account of the history of this problem is given as well.
A similar formula holds in fact also for much larger \(h\) (the best result in this respect is due to Meurman \cite{Meu01}).
For \( k = 3 \), the author \cite{Top16} proved that
\begin{align} \label{eqn: result for D_3(x, h)}
  D_3^\pm(x, h) = x P_{3, h}(\log x) + \BigO{ x^{\frac89 + \varepsilon} } \quad \text{for} \quad h \ll x^\frac23,
\end{align}
with \( P_{3, h} \) a cubic polynomial depending on \(h\).
Both these results rely heavily on spectral methods coming from the theory of automorphic forms in order to obtain the stated error terms.

Unfortunately, the methods used to prove \eqref{eqn: result for D_2(x, h)} and \eqref{eqn: result for D_3(x, h)} do not extend to \( D_k^\pm(x, h) \) with \( k \geq 4 \), the cases we want to focus on in this paper.
Following work by Linnik \cite{Lin63} and Motohashi \cite{Mot80}, Fouvry and Tenenbaum \cite{FT85} used the dispersion method to show that for each \( k \geq 4 \) there exists \( \delta_k > 0 \), such that
\begin{align} \label{eqn: result for D_k(x, h) by Fouvry and Tenenbaum}
  D_k^+(x, 1) = x P_{k, 1}(\log x) + \BigO{ x e^{ -\delta_k \sqrt{\log x} } },
\end{align}
with \( P_{k, 1} \) a polynomial of degree \(k\) (see \cite{Hea86} for more on the history of this problem).
In a recent preprint, Drappeau \cite{Dra15} refined their method and used spectral methods to get a power-saving in the error term -- namely, he showed that there exists \( \delta > 0 \), such that
\begin{align} \label{eqn: result for D_k(x, h) by Drappeau}
  D_k^+(x, h) = x P_{k, h}(\log x) + \BigO{ x^{1 - \frac\delta k} } \quad \text{for} \quad h \ll x^\delta,
\end{align}
with \( P_{k, h} \) a polynomial of order \(k\) depending on \(h\).

Finally, we want to mention the work of Bykovski{\u\i} and Vinogradov \cite{BV87}, where they state a result which is considerably better than \eqref{eqn: result for D_k(x, h) by Drappeau}.
Unfortunately, their proposed proof is incomplete and does not seem to yield the error terms claimed in their paper\footnote{In particular, the step from (5.6) to (5.7) is not correct unless \(n_1\) and \(n_2\) are coprime, and it is unclear how their proposed treatment of \( S(n_1, n_2) \) should work for general \(n_1\) and \(n_2\). See also the comments after Lemma \ref{lemma: main lemma} for another problematic issue.}.
Nevertheless, their initial approach turned out to be useful and led us, together with some new ingredients, to a proof of the following theorem, which improves on \eqref{eqn: result for D_k(x, h) by Fouvry and Tenenbaum} and \eqref{eqn: result for D_k(x, h) by Drappeau}.

\begin{thm} \label{thm: main theorem}
  Let \( k \geq 4 \), and let \( h \) be a positive integer such that \( h \ll x^\frac{15}{19} \).
  Then,
  \[ D_k^\pm(x, h) = x P_{k, h}(\log x) + \BigO{ x^{ 1 - \frac4{15k - 9} + \varepsilon } + x^{ \frac{56}{57} + \varepsilon } }, \]
  where \( P_{k, h} \) is a polynomial of degree \(k\) depending on \(h\), and where the implicit constants depend on \(k\) and \(\varepsilon\).
\end{thm}

The analogous result for the sum weighted by a smooth function is as follows.

\begin{thm} \label{thm: main theorem, smooth version}
  Let \( w : [1/2, 1] \rightarrow [0, \infty) \) be smooth and compactly supported.
  Let \(h\) be a non-zero integer such that \( h \ll x^\frac{15}{19} \).
  Then
  \[ \sum_n w\left( \frac nx \right) d_k(n) d(n \pm h) = x P_{k, h, w}(\log x) + \BigO{ x^{ 1 - \frac1{3k - 2} + \varepsilon} + x^{ \frac{37}{38} + \frac\theta{19} + \varepsilon} }, \]
  where \( P_{k, h, w} \) is a polynomial of degree \(k\) depending on \(w\) and \(h\).
  The implicit constants depend on \(w\), \(k\) and \(\varepsilon\).
\end{thm}

By \(\theta\) we denote the bound in the Ramanujan-Petersson conjecture (see section \ref{subsection: The Kuznetsov formula and large sieve inequalities} for a precise definition -- in any case \( \theta = \frac7{64} \) is admissible).
Both our results can be shown to hold also for larger \(h\), namely up to \( h \ll x^{1 - \varepsilon} \), but for the cost of weaker error terms and with slight changes in the main terms due to the influence of the large shift parameter \(h\).
It seems possible to go even beyond this range, although we have not tried to do so here.

Before we set about to prove these results, we want to describe first in broad terms the approach followed in this paper.
The most direct way to handle shifted convolutions like \( D_k^\pm(x, h) \) is to open one of the divisor functions, and then try to evaluate the arising divisor sums over arithmetic progressions in some way.
In contrast to the mentioned works \cite{Dra15}, \cite{FT85}, \cite{Lin63} and \cite{Mot80}, we have chosen to open \( d_k(n) \) -- although this approach is more difficult from a combinatorial point of view as we have to deal with more variables, the main advantage is that it is much easier to handle \( d(n) \) over arithmetic progressions than \( d_k(n) \) with \( k \geq 4 \).

This way we arrive at sums of the form
\begin{align} \label{eqn: sum after splitting d_k(n)}
  \sum_\twoln{a_1, \ldots, a_k}{a_i \asymp A_i} d(a_1 \cdots a_k + h),
\end{align}
where we can assume that the variables \( a_1, \ldots a_k \) are supported in dyadic intervals \( a_i \asymp A_i \).
The most difficult case occurs when all the variables have small support, which in our case means \( A_i \ll x^{\frac14 + \varepsilon} \), since it is then not enough to just average over one or two of the variables to get an asymptotic formula.
Here we follow the idea of Bykovski{\u\i} and Vinogradov \cite{BV87} and insert the expected main term \( \Phi_0(b) \) for the sum
\[ \Phi(b) := \sum_{a_1 \asymp A_1} d(a_1 b + h) \]
manually into \eqref{eqn: sum after splitting d_k(n)}, so that the latter can be written as
\[ \sum_\twoln{a_2, \ldots, a_k}{a_i \asymp A_i} \Phi_0(a_2 \cdots a_k) - \sum_\twoln{a_2, \ldots, a_k}{a_i \asymp A_i} \left( \Phi_0(a_2 \cdots a_k) - \Phi(a_2 \cdots a_k) \right). \]

While the first sum will be part of the eventual main term, we need to find an upper bound for the second sum.
To do so, we use the Cauchy-Schwarz inequality to bound it by
\[ \left( \sum_{ b \asymp A_2 \cdots A_k } d_{k - 1}(b)^2 \right)^\frac12 \left( \sum_{ b \asymp A_2 \cdots A_k } \left( \Phi_0(b) - \Phi(b) \right)^2 \right)^\frac12, \]
which has the important effect that the variables \( a_2, \ldots, a_k \) are now merged into one large variable \(b\).
After opening the square in the right factor, we are faced with three different sums, the most difficult of them being
\[ \sum_{ b \asymp A_2 \cdots A_k } \Phi(b)^2 = \sum_{a_1, \tilde a_1 \asymp A_1} \sum_{ b \asymp A_2 \cdots A_k } d(a_1 b + h) d(\tilde a_1 b + h). \]

At the heart of our method lies the evaluation of the inner sum on the right hand side, which can be seen as a variation of the binary additive divisor problem and which has been treated by the author in a separate paper \cite{Top15}.
We just want to mention the two main difficulties involved here:
First, we need to evaluate this sum for very large \(a_1\) and \(\tilde a_1\) , namely at least up to the size of \( x^{\frac14 + \varepsilon} \), and the only way to do this with current techniques seems to be by using spectral methods.
Second, while for coprime \(a_1\) and \(\tilde a_1\) this can be done in a fairly standard way, the problem gets much more difficult when \( ( a_1, \tilde a_1 ) > 1 \), and it is a non-trivial issue to employ spectral methods in this situation.

\subsection*{Acknowledgements}

I would like to thank my advisor Valentin Blomer for his support and his valuable advice.
I would also like to thank Sary Drappeau for helpful comments.

\section{Proof of Theorems \ref{thm: main theorem} and \ref{thm: main theorem, smooth version}}

Let \( w : \mathbb{R} \rightarrow [0, \infty) \) be a smooth function compactly supported in \( \left[ 1/2, 1 \right] \) and satisfying
\[ w^{ (\nu) } \ll \frac1{ \Omega^\nu } \quad \text{for} \quad \nu \geq 0, \quad \text{and} \quad \int \! \big| w^{ (\nu) }(\xi) \big| d\xi \ll \frac1{ \Omega^{\nu - 1} } \quad \text{for} \quad \nu \geq 1, \]
for some \( \Omega < 1 \).
We will look at the sum
\[ \Psi(x) := \sum_n w\left( \frac nx \right) d_k(n) d(n + h), \quad h \in \mathbb{Z}, \, h \neq 0, \]
and we will prove an asymptotic formula of the form
\begin{align} \label{eqn: main result for Psi(x)}
  \Psi(x) = M(x) + R(x),
\end{align}
where the main term \( M(x) \) is given by
\[ M(x) := \int \! w\left( \frac\xi x \right) P_{k, h}(\log x, \log \xi,  \log(\xi + h)) \, d\xi, \]
with a polynomial \( P_{k, h} \) of degree \( k \) depending on \(h\), and where we have the following estimate for the error term \( R(x) \), valid for \( h \ll x^{1 - \varepsilon} \),
\begin{align} \label{eqn: main result for R(x)}
  R(x) \ll \frac{ x^{1 - \frac1{3k - 2} + \varepsilon} }{ \Omega^{\frac14 + \frac1{ 4 (3k - 2) } } } + x^{ \frac{37}{38} + \varepsilon } \left( \frac1{\Omega^\frac12} + x^\frac\theta{19} \right) \left( 1 + \left( \frac{ |h| }{ x^\frac{15}{19} } \right)^\frac1{16} \right). 
\end{align}
Remember that \(\theta\) denotes the bound in the Ramanujan-Petersson conjecture (see section \ref{subsection: The Kuznetsov formula and large sieve inequalities}).

The choice \( \Omega = 1 \) directly proves Theorem \ref{thm: main theorem, smooth version}.
Moreover, with the choice \( \Omega = x^{ -\frac1{57} } \) we get
\[ R(x) \ll x^{ \frac{229}{228} - \frac{227}{ 228 (3k - 2) } + \varepsilon } + x^{ \frac{56}{57} + \varepsilon } \left( 1 + \left( \frac{ |h| }{ x^\frac{15}{19} } \right)^\frac1{16} \right), \]
while the choice \( \Omega = x^{ -\frac4{15k - 9} } \) leads to
\[ R(x) \ll x^{ 1 - \frac4{ 15k - 9 } + \varepsilon } + x^{ \frac{37}{38} + \varepsilon } \left( x^\frac2{15k - 9} + x^\frac\theta{19} \right) \left( 1 + \left( \frac{ |h| }{ x^\frac{15}{19} } \right)^\frac1{16} \right). \]
We use the former bound for \( k \leq 15 \) and the latter bound for \( k \geq 16 \), so that 
\[ R(x) \ll x^{ 1 - \frac4{15k - 9} + \varepsilon } + x^{ \frac{56}{57} + \varepsilon } \left( 1 + \left( \frac{ |h| }{ x^\frac{15}{19} } \right)^\frac1{16} \right). \]
After choosing appropriate weight functions, this proves Theorem \ref{thm: main theorem}.

\subsection{Splitting the divisor function}

In order to prove \eqref{eqn: main result for Psi(x)}, we will open \( d_k(n) \) and then dyadically split the supports of the appearing variables.
This will be carried out rigorously in the following subsection -- for the moment, just assume that we have a sum of the form
\begin{align} \label{eqn: definition of splitted Psi}
  \Psi_{v_1, \ldots, v_k} := \sum_{a_1, \ldots, a_k} w\left( \frac{ a_1 \cdots a_k }x \right) v_1(a_1) \cdots v_k(a_k) d( a_1 \cdots a_k + h),
\end{align}
where \( v_1, \ldots, v_k \) are smooth and compactly supported functions satisfying
\[ \supp v_i \asymp A_i \quad \text{and} \quad v_i^{ (\nu) } \ll \frac1{ {A_i}^\nu } \quad \text{for} \quad \nu \geq 0. \]
Our main results are three asymptotic estimates for \( \Psi_{v_1, \ldots, v_k} \), which we state together in the following lemma.

\begin{lemma} \label{lemma: main lemma}
  We have the asymptotic formula
  \[ \Psi_{v_1, \ldots, v_k} = M_{v_1} + R_{v_1}, \]
  where \( M_{v_1} \) is a main term given by
  \[ M_{v_1} := \int \! w\left( \frac\xi x \right) \sum_\twoln{ a_2, \ldots, a_k }{ d \mid a_2 \cdots a_k } \frac{ v_2(a_2) \cdots v_k(a_k) }{a_2 \cdots a_k} v_1\left( \frac\xi{ a_2 \cdots a_k } \right) \lambda_{h, d}(\xi + h) \, d\xi, \]
  with
  \[ \lambda_{h, d}(\xi) := \frac{ c_d(h) }d ( \log\xi + 2\gamma - 2\log d ), \]
  and where we have the following bounds for the error term \( R_{v_1} \), valid for \( h \ll x^{1 - \varepsilon} \),
  \begin{align}
    R_{v_1} &\ll \frac{ x^{\frac32 + \varepsilon} }{ {A_1}^\frac32 \Omega^\frac12 }, \label{eqn: first estimate in main lemma} \\
    R_{v_1} &\ll \frac{ x^{\frac32 + \varepsilon} }{A_1 A_2} \left( \frac1{ \Omega^\frac12 } + \frac{ (A_1 A_2)^{2 \theta} }{ x^{\theta} } \right) \left( 1 + \frac{ {A_2}^\frac12 }{ {A_1}^\frac12 } \right) \left( 1 + |h|^\frac14 \frac{ {A_1}^\frac12 {A_2}^\frac12 }{ x^\frac12 } \right), \label{eqn: second estimate in main lemma} \\
    R_{v_1} &\ll \frac{ x^{1 + \varepsilon} }{ {A_1}^\frac12 \Omega^\frac14} + {A_1}^\frac38 x^{\frac78 + \varepsilon} \left( \frac1{ \Omega^\frac18 } + \frac{ x^\frac\theta4 }{ {A_1}^{\frac34\theta} } + \frac{ |h|^\frac1{16} }{ {A_1}^\frac3{16} } \frac{ x^\frac\theta4 }{ {A_1}^{ \frac34\theta } } \right). \label{eqn: third estimate in main lemma}
  \end{align}
  The implicit constants depend only on \(k\), the involved functions \( w, v_1, \ldots, v_k \) and \(\varepsilon\).
\end{lemma}

When \(A_1\) is large enough that it makes sense to average over \(a_1\) alone, we get the first bound \eqref{eqn: first estimate in main lemma}, which is proven in Section \ref{section: first estimate in main lemma}.
The proof essentially boils down to the evaluation of the following sums over arithmetic progressions modulo \( b = a_2 \cdots a_k \),
\[ \sum_{a_1} w\left( \frac{a_1 b}x \right) v_1(a_1) d(a_1 b + h), \]
for which we can get a non-trivial asymptotic formula as long as \( b \ll x^{\frac23 - \varepsilon} \).
Consequently, also the bound \eqref{eqn: first estimate in main lemma} is non-trivial only for \( A_1 \gg x^{\frac13 + \varepsilon} \).

A further gain in the error term can be achieved here by averaging over another variable \(a_2\), as shown in Section \ref{section: second estimate in main lemma}.
The main ingredient is the Kuznetsov formula that enables us to exploit the cancellation between the Kloosterman sums that arise when the Voronoi summation formula is used to evaluate the sums above.
The resulting bound \eqref{eqn: second estimate in main lemma} is useful when \( A_1 A_2 \gg x^{\frac12 + \varepsilon} \).

The most difficult case occurs when none of the \(A_i\) is particularly large.
It is handled in Section \ref{section: third estimate in main lemma}, and the path we follow there is in some sense dual to the proof of the first bound:
Instead of averaging over \(a_1\), we use the Cauchy-Schwarz inequality to merge the variables \( a_2, \ldots, a_k \) to one large variable \(b\), so that we can then evaluate the sum over this new variable asymptotically.
As mentioned in the introduction, the main difficulty lies in the treatment of the sums
\[ \sum_b w\left( \frac{a_1 b}x \right) w\left( \frac{\tilde a_1 b}x \right) d( a_1 b + h ) d(\tilde a_1 b + h), \]
where \(a_1\) and \( \tilde a_1 \) are of the size \( a_1, \tilde a_1 \asymp A_1 \).
The asymptotic formula we have for these sums has, at best, a non-trivial error term as long as \( a_1, \tilde a_1 \ll x^{\frac13 - \varepsilon} \), and thus the resulting bound \eqref{eqn: third estimate in main lemma} is also non-trivial only for \( A_1 \ll x^{\frac13 - \varepsilon} \).
Note furthermore that this bound is useful only when \( A_i \gg x^\varepsilon \).

Of course, the statement of Lemma \ref{lemma: main lemma} is symmetric in all the variables.
The optimal strategy for given \( A_1, \ldots, A_k \) would be to pick the \(A_i\), which is the largest and which is always at least as large as \( x^\frac1k \), and then apply one of the bounds \eqref{eqn: first estimate in main lemma} or \eqref{eqn: third estimate in main lemma}.
This is essentially the path that Bykovski{\u\i} and Vinogradov \cite{BV87} wanted to take.
Unfortunately, this strategy does not go through, as there is a gap at \( A_i \asymp x^\frac13 \) where both methods fail to give a non-trivial result -- in fact, in the worst case, if for example \( A_1 = A_2 = A_3 \asymp x^\frac13 \) and \( A_4 = \ldots = A_k \asymp 1 \), there is no way to get a non-trivial result from these two bounds alone.

However, we still have another bound at our disposal.
So, if there are two \( A_{i_1}, A_{i_2} \gg x^\frac1k \), at least one of the estimates \eqref{eqn: second estimate in main lemma} or \eqref{eqn: third estimate in main lemma} will always be sufficiently good to get a power saving at the end.
If there is only one \( A_i \gg x^\frac1k \), we can bridge the gap at \( A_i \asymp x^\frac13 \) by using the bound \eqref{eqn: third estimate in main lemma} with respect to one of the other \( A_i \).
More specifically, set
\[ X_1 := \frac{ x^\frac k{3k - 2} }{ \Omega^\frac{k - 1}{ 2 (3k - 2) } }, \quad X_2 := x^\frac5{19} \quad \text{and} \quad X_3 := \left( \frac x{X_1} \right)^\frac1{k - 1} = x^\frac2{3k - 2} \Omega^\frac1{ 2 (3k - 2) }. \]
If one of the \( A_i \) is large enough so that \( A_i \gg X_1 \), we use \eqref{eqn: first estimate in main lemma} to get the estimate
\[ R_{v_1} \ll \frac{ x^{ 1 - \frac1{3k - 2} + \varepsilon } }{ \Omega^{ \frac14 + \frac1{ 4 ( 3k - 2 ) } } }. \]
If there are two \( A_{i_1} \), \( A_{i_2} \) satisfying \( A_{i_1}, A_{i_2} \gg X_2 \), we make use of \eqref{eqn: second estimate in main lemma} and get
\[ R_{v_1} \ll x^{ \frac{37}{38} + \varepsilon } \left( \frac1{ \Omega^\frac12 } + x^\frac\theta{19} \right) \left( 1 + \left( \frac{ |h| }{ x^\frac{18}{19} } \right)^\frac14 \right). \]
Otherwise, there has to be at least one \( A_i \) such that \( X_3 \ll A_i \ll X_2 \), which means that we can use \eqref{eqn: third estimate in main lemma}, hence getting the bound
\[ R_{v_1} \ll \frac{ x^{1 - \frac1{3k - 2} + \varepsilon} }{ \Omega^{\frac14 + \frac1{ 4 (3k - 2) } } } + \frac{ x^{ \frac{37}{38} + \varepsilon } }{ \Omega^\frac18 } + x^{ \frac{37}{38} + \frac\theta{19} + \varepsilon } \left( 1 + \left( \frac{ |h| }{ x^\frac{15}{19} } \right)^\frac1{16} \right). \]
All in all this leads to \eqref{eqn: main result for R(x)}.

\subsection{The main term}

We first describe how to split up the \(k\)-th divisor function so that we can conveniently evaluate the main term at the end.
Let \( u_0 : \mathbb{R} \rightarrow [0, \infty) \) be a smooth and compactly supported function such that
\[ \supp u_0 \subset \left[ \frac12, 2 \right] \quad \text{and} \quad \sum_{ \ell \in \mathbb{Z} } u_\ell(\xi) \equiv 1, \quad \text{where} \quad  u_\ell(\xi) := u_0\left( \frac\xi{ 2^\ell } \right). \]
We set
\[ h_\ell(\xi) := u_\ell \left( \frac\xi{X_3} \right) \quad \text{for} \quad \ell \geq 1, \quad \text{and} \quad h_0(\xi) := \sum_{\ell \leq 0} u_\ell \left( \frac\xi{X_3} \right), \]
and define the sums
\[ \Psi^{ (j) }(x) := \sum_{ a_1, \ldots, a_k } w\left( \frac{ a_1 \cdots a_k }x \right) h_{j_1}(a_1) \cdots h_{j_k}(a_k) d( a_1 \cdots a_k + h ), \]
where \( j = (j_1, \ldots, j_k) \) is a \(k\)-tuple with elements in \( \mathbb{N} = \{ 0, 1, 2, 3, \ldots \} \), so that our main sum can be split up as
\[ \Psi(x) = \sum_\twoln{ j \in \mathbb{N}^k }{ j \neq (0, \ldots, 0) } \Psi^{ (j) }(x). \]

Given a \(k\)-tuple \(j\), there is at least one coordinate \( j_i > 0 \), so that we can use Lemma \ref{lemma: main lemma} with respect to the corresponding variable \( a_i \).
As it turns out, it does not matter which one we choose -- but for the moment we will assume that we can take \( j_1 \).
We split dyadically all the occurring \( h_0(\xi) \) in \( \Psi^{ (j) }(x) \), apply Lemma \ref{lemma: main lemma}, and then sum everything up again, so that
\[ \Psi^{ (j) }(x) = M^{ (j) }(x) + R^{ (j) }(x), \]
where \( R^{ (j) }(x) \) is bounded by \eqref{eqn: main result for R(x)}, and where
\[ M^{ (j) }(x) := \int \! w\left( \frac\xi x \right) \sum_d \frac{ c_d(h) }{ d^2 } \lambda_{h, d}(\xi + h) \sum_\twoln{ j \in \mathbb{N}^k }{ j \neq (0, \ldots, 0) } M_{\text a}^{ (j) }(\xi) \, d\xi, \]
with
\[ M_{\text a}^{ (j) }(\xi) := d \sum_\twoln{ a_2, \ldots, a_k }{ d \mid a_2 \cdots a_k } h_{j_1}\left( \frac\xi{ a_2 \ldots a_k } \right) \frac{ h_{j_2}(a_2) \cdots h_{j_k}(a_k) }{ a_2 \cdots a_k }. \]

We use Mellin inversion to write this sum as
\[ M_{\text a}^{ (j) }(\xi) = \frac1{2\pi i} \int_{ (\sigma) } \! \hat h_{j_1}(s) \frac{d^s}{\xi^s} Z(s; d) \, ds, \quad \sigma < 0, \]
with
\[ \hat h_{j_1}(s) := \int_0^\infty \! h_{j_1}(\eta) \eta^{s - 1} \, d\eta \quad \text{and} \quad Z(s; d) := d^{1 - s} \sum_\twoln{ a_2, \ldots, a_k }{ d \mid a_2 \cdots a_k } \frac{ h_{j_2}(a_2) \cdots h_{j_k}(a_k) }{ ( a_2 \cdots a_k )^{1 - s} }. \]
For integers \( d_2, \ldots, d_k \) such that \( d_2 \cdots d_k = d \) define
\[ c_2 := \frac d{d_2}, \quad c_3 := \frac d{d_2 d_3}, \quad \ldots, \quad c_{k - 1} := \frac d{ d_2 \cdots d_{k - 1} }, \quad c_k := 1. \]
Then we have
\[ \sum_\twoln{a_2, \ldots, a_k}{ d \mid a_2 \cdots a_k } (\ldots) = \sum_{ d_2 \cdots d_k = d } \sum_{ (a_2, d) = d_2 } \sum_{ (a_3, c_2) = d_3 } \cdots \sum_{ ( a_{k - 1}, c_{k - 2} ) = d_{k - 1} } \sum_{ d_k \mid a_k } (\ldots), \]
which means that
\[ Z(s, d) = \sum_{ d_2 \cdots d_k = d } \prod_{i = 2}^k \left( \sum_{ (a, c_i) = 1 } \frac{ h_{j_i}(d_i a) }{ a^{1 - s} } \right). \]

The sums running over \(a\) can be evaluated in the usual way using Mellin inversion and the residue theorem, leading to
\[ \sum_{ (a, c_i) = 1 } \frac{ h_{j_i}(d_i a) }{ a^{1 - s} } = \frac{ \psi_0(c_i) }{ {d_i}^s } H_{j_i}(s) + \BigO{ \frac{ {d_i}^{1 - \varepsilon} }{ \left( 2^{j_i} X_3 \right)^{1 - \varepsilon} } }, \]
where the functions \( H_{j_i}(s) \) are defined as
\begin{align*}
  H_{j_i}(s) &:= \left( 2^{j_i} X_3 \right)^s \int_\frac12^2 \! u_0(\eta) \eta^{s - 1} \, d\eta, \quad \text{for} \quad j_i \geq 1, \\
  H_0(s) &:= \zeta(1 - s) {d_i}^s \frac{ \psi_{-s}(c_i) }{ \psi_0(c_i) } - \frac{ {X_3}^s }s \int_1^2 \! u_0'(\eta) \eta^s \, d\eta,
\end{align*}
with
\[ \psi_\alpha(n) := \prod_{p \mid n} \left( 1 - \frac1{ p^{1 + \alpha} } \right). \]
Because of
\[ \hat h_{j_1}(s) = H_{j_1}(s), \]
we can write \( M_{\text a}^{ (j) }(\xi) \) as
\[ M_{\text a}^{ (j) }(\xi) = \frac1{2\pi i} \sum_{d_2 \cdots d_k = d} \psi_0(c_2) \cdots \psi_0(c_k) \int_{ (\sigma) } \! \prod_{i = 1}^k \left( H_{j_i}(s) + \BigO{ \frac{ {d_i}^{1 - \varepsilon} }{ \left( 2^{j_i} X_3 \right)^{1 - \varepsilon} } } \right) \, \frac{ds}{\xi^s}. \]
Note that this expression is independent of the variable chosen with respect to Lemma \ref{lemma: main lemma}.

At this point we sum all the functions \( H_{j_i}(s) \) with \( j_i \geq 1 \) together, so that
\[ \sum_\twoln{ j \in \mathbb{N}^k }{ j \neq (0, \ldots, 0) } M_{\text a}^{ (j) }(\xi) = \sum_\twoln{ j \in \{ 0, 1 \}^k }{ j \neq (0, \ldots, 0) } M_{\text b}^{ (j) }(\xi) + \BigO{ \frac{ d^{1 - \varepsilon} }{ {X_3}^{1 - \varepsilon} } }, \]
where
\[ M_{\text b}^{ (j) }(\xi) := \frac1{2\pi i} \sum_{ d_2 \cdots d_k = d } \psi_0(c_2) \cdots \psi_0(c_k) \int_{ (\sigma) } \! \frac1{\xi^s} \left( \prod_{i = 1}^k G_{j_i}(s) \right) \, ds, \]
with
\[ G_1(s) := \frac{ {X_3}^s }s \int_1^2 \! u_0'(\eta) \eta^s \, d\eta \quad \text{and} \quad G_0(s) := H_0(s). \]

Next we move the line of integration to \( \sigma = 1 - \varepsilon \), and use the residue theorem to extract a main term from the pole at \( s = 0 \).
Because of
\[ G_1(s) \ll \frac{ {X_3}^{ \Re(s) } }{ |s|^\nu } \quad \text{for} \quad \nu \geq 0, \quad \text{and} \quad \zeta( \varepsilon + it ) \ll |t|^\frac12, \]
we get that
\[ M_{\text b}^{ (j) }(\xi) = P_{k - 1}(\log x, \log \xi) + R_{\text b}^{ (j) }(\xi) + \BigO{ \frac{ {X_3}^{k - 1} }{ x^{1 - \varepsilon} } }, \]
where \( P_{k - 1} \) is a polynomial of degree \( k - 1 \) depending on \(d\) and on \(h\), and where
\[ R_{\text b}^{ (j) }(\xi) = (-1)^{ k - j_1 - \ldots - j_k } \frac1{2\pi i} \sum_{ d_2 \cdots d_k = d } \psi_0(c_2) \cdots \psi_0(c_k) \int_{ (1 - \varepsilon) } \! \frac{ G_1(s)^k }{\xi^s} \, ds. \]
However, because we can assume that
\[ \frac{ {X_3}^k }x \leq \frac14, \]
and because we can move the line of integration to the right as far as we want, we have in fact
\[ R_{\text b}^{ (j) }(\xi) \ll \frac{ d^\varepsilon }x. \]

All in all, the main term of \( \Psi(x) \) is given by
\begin{multline*}
  \sum_\twoln{ j \in \mathbb{N}^k }{ j \neq (0, \ldots, 0) } M^{ (j) }(x) = \int \! w\left( \frac\xi x \right) P_{k, h}( \log x, \log \xi, \log(\xi + h) ) \, d\xi \\
    + \BigO{ \frac{ x^{1 + \varepsilon} }{X_3} + x^\varepsilon {X_3}^{k - 1} },
\end{multline*}
where \( P_{k, h} \) is a polynomial of degree \( k \) depending on \(h\).
Since the error term here is smaller than \eqref{eqn: main result for R(x)}, this proves the asymptotic evaluation claimed in \eqref{eqn: main result for Psi(x)}.

\section{Proof of \eqref{eqn: first estimate in main lemma}} \label{section: first estimate in main lemma}

We write \eqref{eqn: definition of splitted Psi} as
\[ \Psi_{v_1, \ldots, v_k} = \sum_{a_2, \ldots, a_k} v_2(a_2) \cdots v_k(a_k) \Phi(a_2 \cdots a_k), \]
where
\begin{align} \label{eqn: definition of Phi_1}
  \Phi(b) := \sum_{ m \equiv h \smod b } d(m) f(m) = \sum_r d(rb + h) f(rb + h),
\end{align}
with
\[ f(\xi) := w\left( \frac{\xi - h}x \right) v_1\left( \frac{\xi - h}b \right). \]
This divisor sum over an arithmetic progression can be treated using the Voronoi summation formula for the divisor function (see e.g. \cite[Theorem 2.1]{Top16}).
In particular, we have
\begin{align}
  \begin{split} \label{eqn: Voronoi summation formula}
    \Phi(b) &= \frac1b \int \! \Delta_\delta(\xi) f(\xi) \sum_{d \mid b} \frac{ c_d(h) }{ d^{1 + \delta} } \, d\xi \\
      &\qquad -\frac{2 \pi}b \sum_{d \mid b} \sum_{m = 1}^\infty d(m) \frac{ S(h, m; d) }d \int \! Y_0\left( \frac{4 \pi}d \sqrt{m \xi} \right) f(\xi) \, d\xi \\
      &\qquad +\frac4b \sum_{d \mid b} \sum_{m = 1}^\infty d(m) \frac{ S(h, -m; d) }d \int \! K_0\left( \frac{4 \pi}d \sqrt{m \xi} \right) f(\xi) \, d\xi,
  \end{split}
\end{align}
where \( \Delta_{\delta}(\xi) \) is the differential operator defined by
\[ \Delta_{\delta}(\xi) := \left( \log\xi + 2\gamma + 2 \frac\partial{\partial \delta} \right) \bigg|_{\delta = 0}. \]

Note that
\[ \supp f \asymp x \quad \text{and} \quad f^{ (\nu) }(\xi) \ll \frac1{ (x \Omega)^\nu } \quad \text{for} \quad \nu \geq 0, \]
and
\[ \int \! \left| f^{ (\nu) }(\xi) \right| \, d\xi \ll \frac1{ (x \Omega)^{\nu - 1} } \quad \text{for} \quad \nu \geq 1. \]
From \eqref{eqn: Voronoi summation formula}, it follows easily using Weil's bound to bound the Kloosterman sums and the recurrence relations for Bessel functions to bound the integral transforms, that
\begin{align} \label{eqn: the divisor function in arithmetic progressions}
  \Phi(b) = \frac1b \int \! w\left( \frac\xi x \right) v_1\left( \frac\xi b \right) \Delta_\delta(\xi + h) \sum_{d \mid b} \frac{ c_d(h) }{ d^{1 + \delta} } \, d\xi + \BigO{ x^\varepsilon \frac{b^\frac12}{\Omega^\frac12} }
\end{align}
(see \cite[Section 2]{Blo08} for a more detailed treatment).
This formula holds uniformly in \( b\), and eventually leads to
\[ \Psi_{v_1, \ldots, v_k} = M_{v_1} + \BigO{ \frac{ x^{\frac32 + \varepsilon} }{ {A_1}^\frac32 \Omega^\frac12 } }, \]
with \( M_{v_1} \) given as in Lemma \ref{lemma: main lemma}.

\section{Proof of \eqref{eqn: second estimate in main lemma}} \label{section: second estimate in main lemma}

We write \eqref{eqn: definition of splitted Psi} as
\[ \Psi_{v_1, \ldots, v_k} = \sum_{a_3, \ldots, a_k} v_3(a_3) \cdots v_k(a_k) \Phi_2(a_3 \cdots a_k), \]
where
\[ \Phi_2(b) := \sum_{a_2} v_2(a_2) \Phi(a_2 b) = \sum_{a_2} \sum_{ m \equiv h \smod{a_2 b} } d(m) f_2(m, a_2), \]
with
\[ f_2(\xi, a_2) := w\left( \frac{\xi - h}x \right) v_2(a_2) v_1\left( \frac{\xi - h}{a_2 b} \right). \]
We will closely follow \cite[Section 3]{Top16}, where this type of sum is treated in more detail and where we will also borrow in large parts the notation.
For the sum over \(b\) we can again make use of the formula \eqref{eqn: Voronoi summation formula}, leading to
\begin{align*}
  \Phi_2(b) &= \frac1b \sum_{a_2} \frac1{a_2} \int \! \Delta_\delta(\xi) f_2(\xi, a_2) \sum_{d \mid a_2 b} \frac{ c_d(h) }{ d^{1 + \delta} } \, d\xi \\
      &\phantom{ {} = } -\frac{2 \pi}b \sum_{a_2} \frac1{a_2} \sum_{d \mid a_2 b} \sum_{m = 1}^\infty d(m) \frac{ S(h, m; d) }d \int \! Y_0\left( \frac{4 \pi}d \sqrt{m \xi} \right) f_2(\xi, a_2) \, d\xi \\
      &\phantom{ {} = } +\frac4b \sum_{a_2} \frac1{a_2} \sum_{d \mid a_2 b} \sum_{m = 1}^\infty d(m) \frac{ S(h, -m; d) }d \int \! K_0\left( \frac{4 \pi}d \sqrt{m \xi} \right) f_2(\xi, a_2) \, d\xi \\
    &=: M_2(b) - 2\pi R_2^+(b) + 4 R_2^-(b).
\end{align*}
The main term of \( \Psi_{v_1, \ldots, v_k} \) is now given by
\[ M_{v_1} = \sum_{a_3, \ldots, a_k} v_3(a_3) \cdots v_k(a_k) M_2(a_3 \cdots a_k), \]
so that it remains to estimate \( R_2^\pm(b) \).

We change the summation over \(a_2\) and \(d\) as follows,
\[ \sum_\twoln{a_2, d}{d \mid a_2 b} (\ldots) = \sum_\twoln{a_2, c, d}{a_2 b = cd} (\ldots) = \sum_c \sum_{ \frac b{ (b, c) } \mid d } (\ldots), \]
so that
\[ R_2^\pm(b) = \sum_c \frac1c \sum_{m = 1}^\infty d(m) \sum_{ \frac b{ (b, c) } \mid d } \frac{ S(h, \pm m; d) }d F_2^\pm(d, m), \]
with
\[ F_2^\pm(d, m) := \frac1d \int \! B^\pm\left( \frac{4 \pi}d \sqrt{m \xi} \right) f_2\left( \xi, \frac{cd}b \right) \, d\xi, \]
and
\[ B^+(\xi) := Y_0(\xi) \quad \text{and} \quad B^-(\xi) := K_0(\xi). \]
A simple exercise using the recurrence relations for Bessel functions shows that we can cut the sum over \(m\) in \( R_2^-(b) \) and \( R_2^+(b) \) at
\[ M_0^- := \frac{ x^{1 + \varepsilon} }{ (c A_1)^2 } \quad \text{and} \quad M_0^+ := \frac{ x^{1 + \varepsilon} }{ (c A_1 \Omega)^2 }. \]
After furthermore dividing the support of \(m\) dyadically into intervals \( [M, 2M] \), it is hence enough to consider the sums
\[ R_{2 \text a}^\pm(M, c) := \sum_{ M < m \leq 2M } d(m) \sum_{ \frac b{ (b, c) } \mid d } \frac{ S(h, \pm m; d) }d F_2^\pm(d, m). \]

\subsection{The Kuznetsov formula and the large sieve inequalities} \label{subsection: The Kuznetsov formula and large sieve inequalities}

We want to treat the sums \( R_{2 \text a}^\pm(M, c) \) by spectral methods, i.e. the Kuznetsov formula, and in order to do so we need to fix the notation first.
For a background on the theory, we refer to \cite{DI82}.

We will work over the Hecke congruence subgroup \( \Gamma := \Gamma_0(q) \).
Denote by \( \mathfrak{M}_\ell(\Gamma) \) the space of holomorphic cusp forms of weight \(\ell\) with respect to \( \Gamma_0(q) \), and by \( \theta_\ell(q) \) its dimension.
Let \( \{ f_{j, \ell} \}_{ j \leq \theta_\ell(q) } \) be an orthonormal basis of \( \mathfrak{M}_\ell(\Gamma) \).
Then the Fourier expansion of \( f_{j, \ell} \) around \(\infty\) is given by
\[ f_{j, \ell}(z) = \sum_{n = 1}^\infty \psi_{j, \ell}(n) e(nz). \]

Let \( L_\text{cusp}^2( \Gamma \backslash \mathbb{H} ) \) be the Hilbert space of cusp forms with respect to \( \Gamma_0(q) \), and let \( \{ u_j \}_{j \geq 1} \) be an orthonormal basis of this space consisting of Maa{\ss} cusp forms with corresponding eigenvalues \( \left( \frac14 - i \kappa_j \right) \left( \frac14 + i \kappa_j \right) \).
We can assume that either \( i \kappa_j > 0 \) or \( \kappa_j \geq 0 \), depending on whether \(u_j\) is exceptional or not.
Note that if \(u_j\) is exceptional, we know by the work of Kim and Sarnak \cite{Kim03} that \( i \kappa_j \leq \theta \) with
\begin{align*}
  \theta = \frac7{64}
\end{align*}
(although it is conjectured that such exceptional \(u_j\) do no exist, and that therefore \( \theta = 0 \) should be admissible).
The Fourier expansion of \( u_j \) around \(\infty\) is given by
\[ u_j(z) = y^\frac12 \sum_{n \neq 0} \rho_j(n) K_{i \kappa_j} (2\pi |n| y) e(nx). \]
Finally, for the Eisenstein series \( E_{\mathfrak c}(z; s) \) (as defined in \cite[(1.13)]{DI82}), we have a Fourier expansion around \( \infty \) given by
\begin{align*}
  E_\mathfrak{c}(z; s) = \delta_{ \mathfrak{c} \infty } y^s &+ \pi^\frac12 \frac{ \Gamma\left( s - \frac12 \right) }{ \Gamma(s) } \varphi_{ \mathfrak{c}, 0 }(s) y^{1 - s} \\
    &+ 2 y^\frac12 \frac{ \pi^s }{ \Gamma(s) } \sum_{n \neq 0} |n|^{s - \frac12} \varphi_{ \mathfrak{c}, n }(s) K_{s - \frac12}(2 \pi |n| y) e(nx).
\end{align*}

The Kuznetsov formula now reads as follows (see \cite[Theorem 1]{DI82}).

\begin{thm}\label{thm: Kuznetsov formula}
  Let \( f : (0, \infty) \rightarrow \mathbb{C} \) be smooth with compact support and let \(h\), \(m\) be two positive integers.
  Then
  \begin{align*}
    \sum_{ d \equiv 0 \smod q } \frac{ S(h, m; d) }d &f\left( 4\pi \frac{ \sqrt{hm} }d \right) = \sum_{j = 1}^\infty \frac{ \overline{ \rho_j }(h) \rho_j(m) }{ \cosh(\pi \kappa_j) } \hat f( \kappa_j ) \\
      &+ \frac1\pi \sum_\mathfrak{c} \int_{-\infty}^\infty \! \left( \frac hm \right)^{-ir} \overline{ \varphi_{ \mathfrak{c}, h } } \left( \frac12 + ir \right) \varphi_{\mathfrak{c}, m} \left( \frac12 + ir \right) \hat f(r) \, dr \\
      &+ \frac1{2 \pi} \sum_\twoln{ k \equiv 0 \smod 2 }{ 1 \leq j \leq \theta_k(q) } \frac{ i^k (k - 1)! }{ \left( 4\pi \sqrt{hm} \right)^{k - 1} } \overline{ \psi_{j, k} }(h) \psi_{j, k}(m) \tilde f(k - 1),
  \intertext{and}
    \sum_{ d \equiv 0 \smod q } \frac{ S(h, -m; d) }d &f\left( 4\pi \frac{ \sqrt{hm} }d \right) = \sum_{j = 1}^\infty \frac{ \rho_j(h) \rho_j(m) }{ \cosh(\pi \kappa_j) } \check f( \kappa_j ) \\
      &+ \frac1\pi \sum_\mathfrak{c} \int_{-\infty}^\infty \! (hm)^{ir} \varphi_{ \mathfrak{c}, h } \left( \frac12 + ir \right) \varphi_{\mathfrak{c}, m} \left( \frac12 + ir \right) \check f(r) \, dr,
  \end{align*}
  where the Bessel transforms are defined by
  \begin{align*}
    \hat f(r) &= \frac\pi{ \sinh(\pi r) } \int_0^\infty \! \frac{ J_{2ir}(\xi) - J_{-2ir}(\xi) }{2i} f(\xi) \, \frac{d\xi}\xi, \\
    \check f(r) &= \frac4\pi \cosh(\pi r) \int_0^\infty \! K_{2ir}(\xi) f(\xi) \, \frac{d\xi}\xi, \\
    \tilde f(\ell) &= \int_0^\infty \! J_\ell(\xi) f(\xi) \, \frac{d\xi}\xi.
  \end{align*}
\end{thm}

We want to use this formula in order to estimate sums of Kloosterman sums, which means that we need to be able to bound the sums appearing on the spectral side.
In this respect, the following lemma will turn out to be useful (see \cite[Lemma 2.9]{Top16}).

\begin{lemma}\label{lemma: estimates for Fourier coefficients}
  Let \( K \geq 1 \) and \( h \geq 1 \).
  Then
  \begin{align*}
    \sum_{ | \kappa_j | \leq K } \frac{ | \rho_j(h) |^2 }{ \cosh(\pi \kappa_j) } &\ll K^2 + (qKh)^\varepsilon (q, h)^\frac12 \frac{ h^\frac12 }q, \\
    \sum_\mathfrak{c} \int_{-K}^K \! \left| \varphi_{ \mathfrak{c}, h } \left( \frac12 + ir \right) \right|^2 dr &\ll K^2 + (qKh)^\varepsilon (q, h)^\frac12 \frac{ h^\frac12 }q, \\
    \sum_\twoln{2 \leq \ell \leq K, \, 2 \mid \ell}{ 1 \leq j \leq \theta_\ell(q) } \frac{ (\ell - 1)! }{ (4\pi h)^{\ell - 1} } \left| \psi_{j, \ell}(h) \right|^2 &\ll K^2 + (qKh)^\varepsilon (q, h)^\frac12 \frac{ h^\frac12 }q,
  \end{align*}
  where the implicit constants depend only on \( \varepsilon \).
\end{lemma}

Another useful result are the following bounds, which are also known as the large sieve inequalities and are proven in \cite[Theorem 2]{DI82}.

\begin{thm}\label{thm: large sieve inequalities}
  Let \( K \geq 1 \) and \( M \geq \frac12 \) be real numbers, and \( a_m \) a sequence of complex numbers.
  Then
  \begin{align*}
    \sum_{ | \kappa_j | \leq K } \frac1{ \cosh(\pi \kappa_j) } \left| \sum_{M < m \leq 2M} a_n \rho_j(m) \right|^2 &\ll \left( K^2 + \frac{ M^{1 + \varepsilon} }q \right) \sum_{M < m \leq 2M} |a_m|^2, \\
    \sum_\mathfrak{c} \int_{-K}^K \! \left| \sum_{M < m \leq 2M} a_m m^{ir} \varphi_{ \mathfrak{c}, m } \left( \frac12 + ir \right) \right|^2 \, dr &\ll  \left( K^2 + \frac{ M^{1 + \varepsilon} }q \right) \sum_{M < m \leq 2M} |a_m|^2, \\
    \sum_\twoln{2 \leq \ell \leq K, \, 2 \mid \ell}{ 1 \leq j \leq \theta_\ell(q) } \frac{ (\ell - 1)! }{ (4\pi)^{\ell - 1} } \left| \sum_{M < m \leq 2M} a_m m^{ -\frac{\ell - 1}2 } \psi_{j, \ell}(m) \right|^2 &\ll \left( K^2 + \frac{ M^{1 + \varepsilon} }q \right) \sum_{M < m \leq 2M} |a_m|^2,
  \end{align*}
  where the implicit constants depend only on \( \varepsilon \).
\end{thm}

Finally, for the exceptional eigenvalues, we have the following lemma (see \cite[Lemma 2.10]{Top16}).

\begin{lemma} \label{lemma: lemma to treat the exceptional eigenvalues}
  Let \( X, q, h \geq 1 \) be such that \( h^\frac12 X \geq q \).
  Then
  \[ \sum_{ \kappa_j \text{ exc.} } | \rho_j(h) |^2 X^{ 4i \kappa_j } \ll (Xh)^\varepsilon \frac{ h^{2\theta} X^{4\theta} }{ q^{4\theta} } (h, q)^\frac12 \left( 1 + \frac{ h^\frac12 }q \right), \]
  where the implicit constant depends only on \( \varepsilon \).
\end{lemma}

\subsection{An estimate for \( R_{2 \text a}^\pm(M, c) \)}

Going back to our sum, we bring it first into the right shape for the Kuznetsov formula.
We define
\[ \tilde F_2^\pm(d, m) := h(m) F_2^\pm\left( 4\pi \frac{ \sqrt{ |h| m } }d, m \right), \]
where \( h(m) \) is a smooth and compactly supported function such that
\[ h(m) \equiv 1 \quad for \quad m \in [M, 2M], \quad \supp h \asymp M \quad \text{and} \quad h^{ (\nu) }(m) \ll \frac1{ M^\nu }. \]
Moreover, we seperate the variable \(m\) by use of Fourier inversion.
Set
\[ G_0(\lambda) := x^\varepsilon c A_1 \min\left\{ M, \frac1\lambda, \frac1{M \lambda^2} \right\}, \]
and
\[ G_\lambda^\pm(d) := \frac1{ G_0(\lambda) } \int \! \tilde F^\pm(d, m) e(-\lambda m) \, dm, \]
so that
\[ R_{2 \text a}^\pm(M, c) := \int \! G_0(\lambda) \sum_{ M < m \leq 2M } d(m) e(m \lambda) \sum_{ \frac b{ (b, c) } \mid d } \frac{ S(h, \pm m; d) }d G_\lambda^\pm\left( 4\pi \frac{ \sqrt{ |h| m } }d \right) \, d\lambda. \]
Concerning the Bessel transforms of \( G_\lambda^\pm(d) \) appearing in the Kuznetsov formula, we have the following bounds when \( M \ll M_0^- \),
\begin{align}
  \hat G_\lambda(ir), \check G_\lambda(ir) &\ll W^{-2r}, &\text{for} \quad 0 \leq r < \frac14, \label{eqn: first bound for Bessel transforms of G} \\
  \hat G_\lambda(r), \check G_\lambda(r), \tilde G_\lambda(r) &\ll \frac{ r^\varepsilon }{ 1 + r^\frac52 }, &\text{for} \quad r \geq 0, \label{eqn: second bound for Bessel transforms of G} \\
  \intertext{and the following when \( M_0^- \ll M \ll M_0^+ \),}
  \hat G_\lambda(ir), \check G_\lambda(ir) &\ll x^{-\nu}, &\text{for} \quad 0 \leq r < \frac14, \label{eqn: third bound for Bessel transforms of G} \\
  \hat G_\lambda(r), \check G_\lambda(r), \tilde G_\lambda(r) &\ll \frac{x^\varepsilon}{ Z^\frac52 } \left( \frac Zr \right)^\nu, &\text{for} \quad r \geq 0, \label{eqn: fourth bound for Bessel transforms of G}
\end{align}
where we have set
\[ W := \sqrt{ |h| M } \frac{c A_1}x \quad \text{and} \quad Z := \sqrt{xM} \frac{c A_1}x. \]
These estimates can be proven the same way as in \cite[Lemma 3.1]{Top16}.

We will look only at \( R_{2 \text a}^+(M, c) \) in the case \( h > 0 \), since the treatment of the other cases can be done analogously.
We use Theorem \ref{thm: Kuznetsov formula} to get
\[ R_{2 \text a}^\pm(M, c) = \int \! G_0(\lambda) \left( \Xi_\text{exc.}(M) + \Xi_1(M) + \frac1\pi \Xi_2(M) + \frac1{2\pi} \Xi_3(M) \right) \, d\lambda, \]
with
\begin{align*}
  \Xi_1(M) &:= \sum_{j = 1}^\infty \hat G_\lambda^+( \kappa_j ) \left( \frac{ \overline{ \rho_j }(h) }{ \sqrt{ \cosh(\pi \kappa_j) } } \right) \Sigma_j^{ (1) }(M) \\
  \Xi_2(M) &:= \sum_\mathfrak{c} \int_{-\infty}^\infty \! \hat G_\lambda^+(r) \left( h^{-ir} \overline{ \varphi_{ \mathfrak{c}, h } } \left( \frac12 + ir \right) \right) \Sigma_{ \mathfrak{c}, r }^{ (2) }(M) \, dr \\
  \Xi_3(M) &:= \sum_\twoln{ \ell \equiv 0 \smod 2 }{ 1 \leq j \leq \theta_\ell\left( \frac b{ (b, c) } \right) } \tilde G_\lambda^+(\ell - 1) \left( i^\frac \ell2 \sqrt{ \frac{ (\ell - 1)! }{ (4\pi)^{\ell - 1} } } h^\frac{\ell - 1}2 \overline{ \psi_{j, \ell} }(h) \right) \Sigma_{j, \ell}^{ (3) }(M),
  \intertext{and}
  \Sigma_j^{ (1) }(M) &:= \frac1{ \sqrt{ \cosh(\pi \kappa_j) } } \sum_{ M < m \leq 2M } d(m) e(m \lambda) \rho_j(m) \\
  \Sigma_{ \mathfrak{c}, r }^{ (2) }(M) &:= \sum_{ M < m \leq 2M } d(m) e(m \lambda) m^{ir} \varphi_{\mathfrak{c}, m} \left( \frac12 + ir \right) \\
  \Sigma_{j, \ell}^{ (3) }(M) &:= i^\frac\ell2 \sqrt{ \frac{ (\ell - 1)! }{ (4\pi)^{\ell - 1} } } \sum_{ M < m \leq 2M } d(m) e(m \lambda) m^\frac{\ell - 1}2 \psi_{j, \ell}(m).
\end{align*}

First assume \( M \ll M_0^- \).
We divide \( \Xi_1(M) \) into three parts as follows,
\[ \Xi_1(M) = \sum_{\kappa_j \leq 1} (\ldots) + \sum_{1 \leq \kappa_j} (\ldots) + \sum_{ \kappa_j \text{ exc.} } (\ldots) =: \Xi_{1 \text a}(M) + \Xi_{1 \text b}(M) + \Xi_{1 \text c}(M). \]
In \( \Xi_{1 \text a}(M) \) we use \eqref{eqn: second bound for Bessel transforms of G}, Cauchy-Schwarz, and then bound the two arising factors with the help of Lemma \ref{lemma: estimates for Fourier coefficients} and Theorem \ref{thm: large sieve inequalities}, so that
\begin{align}
  \Xi_{1 \text a}(M) &\ll x^\varepsilon \left( \sum_{\kappa_j \leq 1} \frac{ | \rho_j(h) |^2 }{ \cosh(\pi \kappa_j) } \right)^\frac12 \left( \sum_{\kappa_j \leq 1} \left| \Sigma_j^{ (1) }(M) \right|^2 \right)^\frac12 \nonumber \\
    &\ll x^\varepsilon ( b, h )^\frac14 \frac{ (c, b)^\frac12 }c \left( \frac{x^\frac12}{A_1} + \frac{x}{ b^\frac12 {A_1}^2 } \right) \left( 1 + \frac{ h^\frac14 }{b^\frac12} \right). \label{eqn: bound for Xi_1a}
\end{align}
Next, we split \( \Xi_{1 \text b}(M) \) into dyadic segments,
\[ \Xi_1(M, K) := \sum_{K < \kappa_j \leq 2K} \hat G_\lambda^+( \kappa_j ) \left( \frac{ \overline{ \rho_j }(h) }{ \sqrt{ \cosh(\pi \kappa_j) } } \right) \Sigma_j^{ (1) }(M), \]
for which we get, in the same way as above, the bound
\[ \Xi_1(M, K) \ll x^\varepsilon \frac{ (c, b)^\frac12 }c \left( \frac{x^\frac12}{A_1} + \frac{x}{ b^\frac12 {A_1}^2 K } \right) \left( 1 + ( b, h )^\frac14 \frac{ h^\frac14 }{b^\frac12 K} \right) \frac1{K^\frac12}. \]
The estimate we get from this for \( \Xi_{1 \text b}(M) \) is the same as \eqref{eqn: bound for Xi_1a}.
Finally, by using the estimate \eqref{eqn: first bound for Bessel transforms of G} and Lemma \ref{lemma: lemma to treat the exceptional eigenvalues}, we get the following bound for \( \Xi_{1 \text c}(M) \),
\begin{align*}
  \Xi_{1 \text c}(M) &\ll \left( \sum_{ \kappa_j \text{ exc.} } | \rho_j (h) |^2 W^{-4 i \kappa_j} \right)^\frac12 \left( \sum_{ \kappa_j \text{ exc.} } \left| \Sigma_j^{ (1) }(M) \right|^2 \right)^\frac12 \\
    &\ll x^\varepsilon (b, h)^\frac14 \frac{ (c, b) }c \frac{ x^\theta }{ b^{2\theta} } \left( \frac{ x^\frac12 }{A_1} + \frac x{b^\frac12 {A_1}^2 } \right) \left( 1 + \frac{ h^\frac14 }{b^\frac12} \right).
\end{align*}
All in all this leads to
\[ \int \! G_0(\lambda) \Xi_1(M) \, d\lambda \ll (b, h)^\frac14 (c, b) \frac{ x^{\theta + \varepsilon} }{ b^{2\theta} } \left( x^\frac12 + \frac{x}{ b^\frac12 A_1 } \right) \left( 1 + \frac{ h^\frac14 }{b^\frac12} \right). \]

In the other case, when \( M_0^- \ll M \ll M_0^+ \), we split \( \Xi_1(M) \) in the following way,
\[ \Xi_1(M) = \sum_{\kappa_j \leq Z} (\ldots) + \sum_{Z \leq \kappa_j} (\ldots) + \sum_{ \kappa_j \text{ exc.} } (\ldots), \]
and then treat these sums in the same manner as shown above, the main differences being that we now have to use the bounds \eqref{eqn: third bound for Bessel transforms of G} and \eqref{eqn: fourth bound for Bessel transforms of G}, and that the exceptional eigenvalues pose no problems here.
As a result we get
\[ \int \! G_0(\lambda) \Xi_1(M) \, d\lambda \ll (b, h)^\frac14 (c, b)^\frac12 \frac{ x^\varepsilon }{ \Omega^\frac12 } \left( x^\frac12 + \frac{x}{ b^\frac12 A_1 } \right) \left( 1 + \frac{ h^\frac14 }{b^\frac12} \right). \]

The sums \( \Xi_2(M) \) and \( \Xi_3(M) \) can be handled in very much the same way, so we will refrain from giving the details.
Eventually we get the bound
\[ R_2^\pm(b) \ll (b, h)^\frac14 x^\varepsilon \left( \frac1{ \Omega^\frac12 } + \frac{x^\theta}{ b^{2 \theta} } \right) \left( x^\frac12 + \frac x{ b^\frac12 A_1 } \right) \left( 1 + \frac{h^\frac14}{ b^\frac12 } \right), \]
which then immediately leads to \eqref{eqn: second estimate in main lemma}.

\section{Proof of \eqref{eqn: third estimate in main lemma}} \label{section: third estimate in main lemma}

We write
\[ \Psi_{v_1, \ldots, v_k} = \sum_b \delta(b) \Phi(b), \]
where \( \Phi(b) \) is defined as in \eqref{eqn: definition of Phi_1}, and where
\[ \delta(b) := \sum_\twoln{a_2, \ldots, a_k}{ a_2 \cdots a_k = b } v_2(a_2) \cdots v_k(a_k). \]
Furthermore set
\[ B := A_2 \cdots A_k \asymp \frac x{A_1}. \]
If \(B\) is too large, it does not make sense to evaluate the divisor sum over arithmetic progressions as in the sections before.
Instead, we insert the main term from \eqref{eqn: Voronoi summation formula}, namely
\[ \Phi_0(b) := \frac1b \int \! \Delta_\delta(\xi + h) w\left( \frac\xi x \right) v_1\left( \frac\xi b \right) \sum_{d \mid b} \frac{ c_d(h) }{ d^{1 + \delta} } \, d\xi, \]
manually in our sum,
\[ \Psi_{v_1, \ldots, v_k} = \sum_b \delta(b) \Phi_0(b) - \sum_b \delta(b) ( \Phi_0(b) - \Phi(b) ). \]
The main term of \( \Psi_{v_1, \ldots, v_k} \) is then given by the left-most sum -- in fact,
\[ \sum_b \delta(b) \Phi_0(b) = M_{v_1}, \]
where \( M_{v_1} \) is defined as in Lemma \ref{lemma: main lemma}.

It remains to show that the remainder
\[ R_{v_1} := \sum_b \delta(b) ( \Phi_0(b) - \Phi(b) ) \]
is small, and as a first step we use Cauchy-Schwarz,
\[ R_{v_1} \leq \left( \sum_{b \asymp B} | \delta(b) |^2 \right)^\frac12 \left( \sum_b \left| \Phi_0(b) - \Phi(b) \right|^2 \right)^\frac12. \]
While the first factor can be estimated trivially,
\[ \sum_{b \asymp B} |\delta(b) |^2 \ll x^\varepsilon B \ll \frac{ x^{1 + \varepsilon} }{A_1}, \]
the other factor needs more work.
We write
\[ \sum_b \left| \Phi_0(b) - \Phi(b) \right|^2 = \Sigma_1 - 2 \Sigma_2 + \Sigma_3, \]
with
\begin{align*}
  \Sigma_1 := \sum_b \Phi_0(b)^2, \quad \Sigma_2 &:= \sum_b \Phi_0(b) \Phi(b), \quad \Sigma_3 := \sum_b \Phi(b)^2.
\end{align*}
In what follows, we will evaluate these sums and show that
\begin{align}
  \Sigma_1 &= M_0 + \BigO{ x^\varepsilon {A_1}^2 }, \label{eqn: evaluation of Sigma_1} \\
  \Sigma_2 &= M_0 + \BigO{ x^{1 + \varepsilon} + \frac{ x^{\frac13 + \varepsilon} {A_1}^2 }{\Omega^\frac13} }, \label{eqn: evaluation of Sigma_2} \\
  \Sigma_3 &= M_0 + \BigO{ \frac{ x^{1 + \varepsilon} }{\Omega^\frac12} + {A_1}^\frac74 x^{\frac34 + \varepsilon} \left( \frac1{ \Omega^\frac14 } + \frac{ x^\frac\theta2 }{ {A_1}^{\frac32\theta} } + \frac{ |h|^\frac18 }{ {A_1}^\frac38 } \frac{ x^\frac\theta2 }{ {A_1}^{ \frac32\theta } } \right) }, \label{eqn: evaluation of Sigma_3}
\end{align}
where \( M_0 \) is defined in \eqref{eqn: definition of the main term of Sigma_i}.
Hence
\[ R_{v_1} \ll \frac{ x^{1 + \varepsilon} }{ {A_1}^\frac12 \Omega^\frac14} + {A_1}^\frac38 x^{\frac78 + \varepsilon} \left( \frac1{ \Omega^\frac18 } + \frac{ x^\frac\theta4 }{ {A_1}^{\frac34\theta} } + \frac{ |h|^\frac1{16} }{ {A_1}^\frac3{16} } \frac{ x^\frac\theta4 }{ {A_1}^{ \frac34\theta } } \right), \]
thus proving \eqref{eqn: third estimate in main lemma}.

\subsection{Evaluation of \( \Sigma_1 \)}

We have
\[ \Sigma_1 = \iint \! \Delta_{\delta_1}(\xi_1 + h) \Delta_{\delta_2}(\xi_2 + h) w\left( \frac{\xi_1}x \right) w\left( \frac{\xi_2}x \right) \Sigma_{ 1 \text{a} }(\xi_1, \xi_2) \, d\xi_1 d\xi_2, \]
where
\[ \Sigma_{ 1 \text{a} }(\xi_1, \xi_2) := \sum_b \frac{ f_1(b) }{b^2} \sum_{d_1, d_2 \mid b} \frac{ c_{d_1}(h) c_{d_2}(h) }{ {d_1}^{1 + \delta_1} {d_2}^{1 + \delta_2} }, \]
with
\[ f_1(\eta) := v_1\left( \frac{\xi_1}\eta \right) v_1\left( \frac{\xi_2}\eta \right). \]
We use Mellin inversion to evaluate the sum over \(b\), so that we can write
\begin{align} \label{eqn: Sigma_1 after Mellin inversion}
  \Sigma_{ 1 \text{a} }(\xi_1, \xi_2) = \frac1{2\pi i} \int_{ (\sigma) } \! \hat f_1(s) Z_1(s) \, ds,
\end{align}
where \( \sigma > - 1 \), and where
\[ \hat f_1(s) := \int_0^\infty \! f_1(\eta) \eta^{s - 1} \, d\eta \quad \text{and} \quad Z_1(s) := \sum_b \frac1{ b^{2 + s} } \sum_{d_1, d_2 \mid b} \frac{ c_{d_1}(h) c_{d_2}(h) }{ {d_1}^{1 + \delta_1} {d_2}^{1 + \delta_2} }. \]
The Dirichlet series \( Z_1(s) \) converges absolutely for \( \Re(s) > - 1 \), but it is not hard to find an analytic continuation up to \( \Re(s) > -2 \), namely
\[ Z_1(s) = \zeta(2 + s) \sum_{d_1, d_2} \frac{ c_{d_1}(h) c_{d_2}(h) (d_1, d_2)^{2 + s} }{ {d_1}^{3 + \delta_1 + s} {d_2}^{3 + \delta_2 + s} }. \]

We want to move the line of integration in \eqref{eqn: Sigma_1 after Mellin inversion} to \( \sigma = -2 + \varepsilon \) and use the residue theorem to extract a main term.
\( Z_1(s) \) has a pole at \( s = -1 \) with residue
\[ \Res{s = -1} \, Z_1(s) = \sum_{d_1, d_2} \frac{ c_{d_1}(h) c_{d_2}(h) (d_1, d_2) }{ {d_1}^{2 + \delta_1} {d_2}^{2 + \delta_2} } =: C_{\delta_1, \delta_2}(h). \]
Furthermore, we have the bound
\[ Z_1(-2 + \varepsilon + it) \ll x^\varepsilon |t|^{\frac12 + \varepsilon}, \]
which also holds for its derivatives with respect to \( \delta_1 \) and \( \delta_2 \), as well as
\[ \hat f_1(s) \ll B^{ \Re(s) } \min\left\{ 1, \frac1{ |s| }, \frac1{ |s| |s + 1| } \right\}. \]
It follows that
\[ \Sigma_{ 1 \text{a} }(\xi_1, \xi_2) = C_{\delta_1, \delta_2}(h) \int \! \frac{ f_1(\eta) }{\eta^2} \, d \eta + \BigO{ \frac{x^\varepsilon}{B^2} }. \]
One can check that \( C_{\delta_1, \delta_2}(h) \) can be written as
\begin{align} \label{eqn: defintion of C(h)}
  C_{\delta_1, \delta_2}(h) = C_{\delta_1, \delta_2} \gamma_{\delta_1, \delta_2}(h)
\end{align}
where we have set \( C_{\delta_1, \delta_2} := C_{\delta_1, \delta_2}(1) \), and where \( \gamma_{\delta_1, \delta_2}(h) \) is a multiplicative function defined on prime powers by
\begin{align*}
  \gamma_{\delta_1, \delta_2}\left( p^k \right) &:= \sum_{i = 0}^k \frac1{ p^{i + i \delta_1 + i \delta_2} } + \sum_{i = 0}^{k - 1} \sum_{j = 0}^i \left( \frac1{ p^{ (i + 1) + j \delta_1 + (i + 1) \delta_2 } } + \frac1{ p^{ (i + 1) + (i + 1) \delta_1 + j \delta_2 } } \right) \\
    &\phantom{ {} := } - \frac{p - 1}{ p^{3 + \delta_1 + \delta_2} - p^{1 + \delta_1} - p^{1 + \delta_2} + 1 } \sum_{i = 0}^{k - 1} \sum_{j = 0}^i \left( \frac1{ p^{i + j \delta_1 + i \delta_2} } + \frac1{ p^{i + i \delta_1 + j \delta_2} } \right).
\end{align*}
Hence
\[ \Sigma_1 = M_0 + \BigO{ x^\varepsilon {A_1}^2 }, \]
where
\begin{align} \label{eqn: definition of the main term of Sigma_i}
  M_0 := \iiint \! \Delta_{\delta_1}(\xi_1 + h) \Delta_{\delta_2}(\xi_2 + h) C_{\delta_1, \delta_2} \gamma_{\delta_1, \delta_2}(h) F(\xi_1, \xi_2, \eta) \, d\xi_1 d\xi_2 d\eta,
\end{align}
with
\begin{align} \label{eqn: definition of F(xi_1, xi_2, eta)}
  F(\xi_1, \xi_2, \eta) := \frac1{\eta^2} v_1\left( \frac{\xi_1}{\eta} \right) v_1\left( \frac{\xi_2}{\eta} \right) w\left( \frac{\xi_1}x \right) w\left( \frac{\xi_2}x \right).
\end{align}

\subsection{Evaluation of \( \Sigma_2 \)}

We have
\[ \Sigma_2 = \int \! \Delta_{\delta_1}(\xi_1 + h) w\left( \frac{\xi_1}x \right) \sum_{d_1 \leq D_1} \frac{ c_{d_1}(h) }{ {d_1}^{2 + \delta_1} } \Sigma_{ 2 \text{a} }(\xi_1; d_1) \, d\xi_1 + \BigO{ \frac{ x^{1 + \varepsilon} A_1 }{D_1} }, \]
where we have cut the sum over \(d_1\) at \( D_1 \), and where
\[ \Sigma_{ 2 \text{a} }(\xi_1; d_1) := d_1 \sum_r r \sum_{ m \equiv h \smod{r d_1} } f_2( m - h, r ) d(m), \]
with
\[ f_2(\xi_2, \eta) := \frac{ v_1(\eta) }{\xi_2} w\left( \frac{\xi_2}x \right) v_1\left( \frac{\eta \xi_1}{\xi_2} \right). \]
Here we can directly use \eqref{eqn: the divisor function in arithmetic progressions}, which leads to
\[ \Sigma_{ 2 \text{a} }(\xi_1; d_1) := \int \! \Delta_{\delta_2}(\xi_2 + h) \Sigma_{ 2 \text{b} }(\xi_2; d_1) \, d\xi_2 + \BigO{ x^\varepsilon \frac{ {d_1}^\frac32 {A_1}^\frac52 }{x \Omega^\frac12} }, \]
with
\[ \Sigma_{ 2 \text{b} }(\xi_2; d_1) := \sum_r f_2(\xi_2, r) \sum_{d_2 \mid r d_1} \frac{ c_{d_2}(h) }{ {d_2}^{1 + \delta_2} }. \]
We can now evaluate the sum over \(r\) using Mellin inversion in the same way as in the section before.
We have for \( \sigma > 1 \),
\[ \Sigma_{ 2 \text{b} }(\xi_2; d_1) = \frac1{2\pi i} \int_{ (\sigma) } \! \hat f_2(\xi_2, s) Z_2(s) \, ds, \]
where
\[ \hat f_2(\xi_2, s) := \int \! f_2(\xi_2, \eta) \eta^{s - 1} \, d\eta, \]
and
\[ Z_2(s) := \sum_r \frac1{r^s} \sum_{d_2 \mid r d_1} \frac{ c_{d_2}(h) }{ {d_2}^{1 + \delta_2} } = \zeta(s) \sum_{d_2} \frac{ c_{d_2}(h) (d_1, d_2)^s }{ {d_2}^{1 + s+ \delta_2} }. \]
After moving the line of integration to \( \sigma = \varepsilon \) and using the residue theorem, we get
\[ \Sigma_{ 2 \text{b} }(\xi_1, \xi_2; d_1) = \sum_{d_2} \frac{ c_{d_2}(h) (d_1, d_2) }{ {d_2}^{2 + \delta_2} } \int \! f_2(\xi_2, \eta) \, d\eta + \BigO{ \frac{x^\varepsilon}{x} }, \]
which then leads to
\begin{multline*}
  \Sigma_{2 \text a}(\xi_1; d_1) = \iint \! \Delta_{\delta_2}(\xi_2 + h) f_2(\xi_2, \eta) \sum_{d_2} \frac{ c_{d_2}(h) (d_1, d_2) }{ {d_2}^{2 + \delta_2} } \, d\eta d\xi_2 \\
    + \BigO{ x^\varepsilon + x^\varepsilon \frac{ {d_1}^\frac32 {A_1}^\frac52 }{ x \Omega^\frac12 } }.
\end{multline*}
We complete the sum over \(d_1\) again, and  get eventually
\[ \Sigma_2 = M_0 + \BigO{ x^{1 + \varepsilon} + \frac{ x^{1 + \varepsilon} A_1 }{D_1} + x^{\varepsilon} \frac{ {D_1}^\frac12 {A_1}^\frac52 }{\Omega^\frac12} }. \]
The optimal value for \( D_1 \) is
\[ D_1 = \frac{ x^\frac23 \Omega^\frac13 }{A_1}, \]
which gives \eqref{eqn: evaluation of Sigma_2}.

\subsection{Evaluation of \( \Sigma_3 \)}

We have
\[ \Sigma_3 = \sum_{r_1, r_2} v_1(r_1) v_1(r_2) \Sigma_{3 \text a}(r_1, r_2), \]
with
\[ \Sigma_{3 \text a}(r_1, r_2) := \sum_b w\left( \frac{r_1 b}x \right) w\left( \frac{r_2 b}x \right) d(r_1 b + h) d(r_2 b + h). \]
The last sum can be seen as a variant of the binary additive divisor problem, and has been treated in the case \( r_1 \neq r_2 \) in \cite{Top15}.
As a result, we know that \( \Sigma_{3 \text a}(r_1, r_2) \) can be written asymptotically as
\begin{align} \label{eqn: asymptotic evaluation for Sigma_3a}
  \Sigma_{3 \text a}(r_1, r_2) = M_{3 \text a}(r_1, r_2) + R_{3 \text a}(r_1, r_2),
\end{align}
with a main term \( M_{3 \text a}(r_1, r_2) \) and an error term \( R_{3 \text a}(r_1, r_2) \).
More precisely, the main term has the form
\[ M_{3 \text a}(r_1, r_2) := \int \! w\left( \frac{r_1 \eta}x \right) w\left( \frac{r_2 \eta}x \right) \Delta_{\delta_1}(r_1 \eta + h) \Delta_{\delta_2}(r_2 \eta + h) C_{3 \text a}(r_1, r_2, h) \, d\eta, \]
where
\begin{multline*}
  C_{3 \text a}(r_1, r_2, h) := \sum_\twoln{ u_1^\ast \mid (r_1, h) }{ u_2^\ast \mid (r_2, h) } \frac{ {u_1^\ast}^{\delta_1} {u_2^\ast}^{\delta_2} }{ (r_1, h)^{\delta_1} (r_2, h)^{\delta_2} } \psi_{\delta_1}\left( \frac{r_1 u_1^\ast}{ (r_1, h) } \right) \psi_{\delta_2}\left( \frac{r_2 u_2^\ast}{ (r_2, h) } \right) \\
    \cdot \sum_{ \left( d, \frac{ r_1 u_1^\ast r_2 u_2^\ast }{ (r_1, h) (r_2, h) } \right) = 1 } \frac1{ d^{2 + \delta_1 + \delta_2} } c_d\left( \frac{ r_1 h - r_2 h }{ (r_1, h) (r_2, h) } \right),
\end{multline*}
with \( \psi_\alpha(n) \) an arithmetic function defined by
\[ \psi_\alpha(n) := \prod_{p \mid n} \left( 1 - \frac1{ p^{1 + \alpha} } \right). \]
Concerning the error term, we know from \cite[(3.4)]{Top15},
\begin{multline} \label{eqn: main bound for R_3a}
  R_{3 \text a}(r_1, r_2) \ll (r_1, r_2)^\star {A_1}^\frac12 x^{\frac12 + \varepsilon} \\
    \cdot \left( \frac1{\Omega^\frac12} + \left( r_1 r_2, h \frac{ r_2 - r_1 }{ (r_1, r_2) } \right)^{\frac14 + \theta} \frac{ x^\theta }{ {A_1}^{3\theta} } \left( 1 + \frac{ |h|^\frac14 }{ {A_1}^\frac34 } \right) \right),
\end{multline}
where
\[ (r_1, r_2)^\star := \min\left\{ \left( r_1, {r_2}^\infty \right), \left( r_2, {r_1}^\infty \right) \right\}. \]
Of course, there is also the trivial bound
\begin{align} \label{eqn: trivial bound for R_3a}
  R_{3 \text a}(r_1, r_2) \ll x^\varepsilon B.
\end{align}

Since the contribution of the diagonal elements \( r_1 = r_2 \) is negligible, we can bound the respective sums trivially.
Otherwise we use the asymptotic formula \eqref{eqn: asymptotic evaluation for Sigma_3a}, so that
\[ \Sigma_3 = M_3 + R_3, \]
where
\[ M_3 := \sum_{r_1 \neq r_2} v_1(r_1) v_1(r_2) M_{3 \text a}(r_1, r_2), \]
and
\[ R_3 \ll \sum_{r \asymp A_1} | \Sigma_{3 \text a}(r, r) | + \sum_\twoln{ r_1, r_2 \asymp A_1, \,\, r_1 \neq r_2 }{ (r_1, r_2)^\star > R_0} | R_{3 \text a}(r_1, r_2) | + \sum_\twoln{ r_1, r_2 \asymp A_1, \,\, r_1 \neq r_2 }{ (r_1, r_2)^\star \leq R_0} | R_{3 \text a}(r_1, r_2) |, \]
with \( R_0 \ll A_1 \) some constant to be determined at the end.
For the first sum we have
\[ \sum_{r \asymp A_1} | \Sigma_{3 \text a}(r, r) | \ll x^{1 + \varepsilon}. \]
For the second sum, we use the trivial bound \eqref{eqn: trivial bound for R_3a},
\[ \sum_\twoln{r_1, r_2 \asymp A_1, \,\, r_1 \neq r_2}{ (r_1, r_2)^\star > R_0} | R_{3 \text a}(r_1, r_2) | \ll x^\varepsilon B \sum_\twoln{r_1, r_2 \asymp A_1}{ (r_1, r_2) > R_0 } 1 \ll x^\varepsilon B \sum_{r_0 > R_0} \frac{ {A_1}^2 }{ {r_0}^2 } \ll x^{1 + \varepsilon} \frac{A_1}{R_0}. \]
Finally, for the third sum we use \eqref{eqn: main bound for R_3a}.
Note hereby, that
\[ \sum_\twoln{r_1, r_2 \asymp A_1}{ (r_1, r_2)^\star \leq R_0} (r_1, r_2)^\star \ll \sum_{r_0 \leq R_0} r_0 \sum_\twoln{r_1, r_2 \asymp A_1}{ \left( r_1, {r_2}^\infty \right) = r_0 } 1\ll \sum_{r_0 \leq R_0} r_0 \sum_\twoln{r_1, r_2 \asymp A_1}{r_0 \mid r_1} 1 \ll R_0 {A_1}^2, \]
and moreover that,
\begin{align*}
  \sum_\thrln{ r_1, r_2 \asymp A_1 }{r_1 \neq r_2}{ (r_1, r_2)^\star \leq R_0 } & (r_1, r_2)^\star \left( r_1 r_2, h \frac{ r_2 - r_1 }{ (r_1, r_2) } \right)^{\frac14 + \theta} \ll \sum_{r_0 \leq R_0} r_0 \sum_\thrln{ r_0 r_1 \asymp A_1 }{ r_2 \asymp A_1 }{ r_0 r_1 \neq r_2 } \left( r_0 r_1 r_2, h \frac{ r_2 - r_0 r_1 }{ (r_0 r_1, r_2) } \right)^\frac12 \\
    &\ll \sum_{r_0 \leq R_0} r_0 \sum_{t \ll A_1}  \sum_\thrln{ \frac{r_1 r_0}{ (r_0, t) } \asymp \frac{ A_1 }t \asymp r_2 }{ \frac{r_0 r_1}{ (r_0, t) } \neq r_2 }{ \left( \frac{r_0 r_1}{ (r_0, t) }, r_2 \right) = 1 } \left( t^2, \left( r_2 - \frac{r_0 r_1}{ (r_0, t) } \right) h \right)^\frac12 \left( \frac{r_0 r_1 r_2}{ (r_0, t) }, h \right)^\frac12 \\
    &\ll \sum_{s \leq R_0} \sum_\thrln{ r_0 \leq \frac{R_0}s }{ t \ll \frac{A_1}s }{ (r_0, t) = 1 } r_0 s \sum_\twoln{ r_1 r_0 \asymp \frac{ A_1 }{st} \asymp r_2 }{ r_0 r_1 \neq r_2 } ( s t, ( r_2 - r_0 r_1 ) h ) ( r_0 r_1 r_2, h )^\frac12 \\
    &\ll \ldots,
\end{align*}
and after dividing the ranges of the variables \(s\) and \(t\) dyadically into ranges \( s \asymp S \) and \( t \asymp T \),
\begin{align*}
  \ldots &\ll \sum_{S, T} S \sum_{ r_0 \ll \frac{R_0}S } r_0 \sum_\twoln{ r_1 r_0 \asymp \frac{ A_1 }{ST} \asymp r_2 }{ r_0 r_1 \neq r_2 } ( r_0 r_1 r_2, h )^\frac12 \sum_\twoln{s \asymp S}{t \asymp T} ( s t, ( r_2 - r_0 r_1 ) h ) \\
    &\ll |h|^\varepsilon {A_1}^\varepsilon \sum_{S, T} S^2 T \sum_{ r_0 \ll \frac{R_0}S } r_0 (r_0, h)^\frac12 \sum_{ r_1 \asymp \frac{ A_1 }{r_0 ST} } ( r_1, h )^\frac12 \sum_{ r_2 \asymp \frac{A_1}{ST} } ( r_2, h )^\frac12 \\
    &\ll |h|^\varepsilon R_0 {A_1}^{2 + \varepsilon}.
\end{align*}
Hence
\[ \sum_\twoln{ r_1, r_2 \asymp A_1, \,\, r_1 \neq r_2 }{ (r_1, r_2)^\star \leq R_0} | R_{3 \text a}(r_1, r_2) | \ll R_0 {A_1}^\frac52 x^{\frac12 + \varepsilon} \left( \frac1{\Omega^\frac12} + \frac{x^\theta}{ {A_1}^{3 \theta} } \left( 1 + \frac{ |h|^\frac14 }{ {A_1}^\frac34 } \right) \right). \]

We set
\[ R_0 = \min\left\{ A_1, \frac{ x^\frac14 }{ {A_1}^\frac34 } \left( \frac1{\Omega^\frac12} + \frac{x^\theta}{ {A_1}^{3 \theta} } \left( 1 + \frac{ |h|^\frac14 }{ {A_1}^\frac34 } \right) \right)^{-\frac12} \right\}, \]
which leads to
\[ \Sigma_3 = M_3 + \BigO{ x^{1 + \varepsilon} + {A_1}^\frac74 x^{\frac34 + \varepsilon} \left( \frac1{ \Omega^\frac14 } + \frac{ x^\frac\theta2 }{ {A_1}^{\frac32 \theta} } + \frac{ |h|^\frac18 }{ {A_1}^\frac38 }  \frac{ x^\frac\theta2 }{ {A_1}^{\frac32 \theta} } \right) }. \]
It remains to evaluate the main term \(M_3\).

\subsection{The main term of \(\Sigma_3\)}

The main term is given by
\begin{multline*}
  M_3 := \int \! \Delta_{\delta_1}(1) \Delta_{\delta_2}(1) \frac1{ h^{ \delta_1 + \delta_2 } } \!\! \sum_\twoln{u_1, u_2 \mid h}{ u_1^\ast \mid \frac h{u_1}, \,\, u_2^\ast \mid \frac h{u_2} } \!\! {u_1}^{\delta_1} {u_2}^{\delta_2} {u_1^\ast}^{\delta_1} {u_2^\ast}^{\delta_2} \\
    \cdot \sum_\twoln{ d \leq D_0 }{ (d, u_1^\ast u_2^\ast) = 1 } \frac{ M_{3 \text b}(d) }{ d^{2 + \delta_1 + \delta_2} } \, d\eta  + \BigO{ \frac{ A_1 x^{1 + \varepsilon} }{D_0} },
\end{multline*}
where we have cut the sum over \(d\) at \( D_0 \ll A_1 \), and where
\[ M_{3 \text b}(d) := \sum_\thrln{r_1, r_2}{ ( r_1, u_1 d ) = 1 }{ ( r_2, u_2 d ) = 1 } f_{3, 1}\left( \frac{h r_1}{u_1} \right) f_{3, 2}\left( \frac{h r_2}{u_2} \right) \psi_{\delta_1}( r_1 u_1^\ast ) \psi_{\delta_2}( r_2 u_2^\ast ) c_d( r_1 u_2 - r_2 u_1 ), \]
with
\[ f_{3, i}(\xi) := \left( \xi \eta + h \right)^\frac{\delta_i}2 v_1(\xi) w\left( \frac{\xi \eta}x \right). \]

We open the Ramanujan sum, so that
\begin{multline*}
  M_{3 \text b}(d) = \sum_\twoln{y \smod d}{ (y, d) = 1 } \sum_{ ( r_1, u_1 d ) = 1 } f_{3, 1}\left( \frac{h r_1}{u_1} \right) \psi_{\delta_1}( r_1 u_1^\ast ) e\left( \frac{y r_1 h_2}{d_2} \right) \\
    \cdot \sum_{ \left( r_2, u_2 d \right) = 1 } f_{3, 2}\left( \frac{h r_2}{u_2} \right) \psi_{\delta_1}( r_2 u_2^\ast ) e\left( -\frac{y r_2 h_1}{d_1} \right),
\end{multline*}
where we have set
\[ d_1 := \frac d{ ( d, u_1 ) }, \quad h_1 := \frac{u_1}{ ( d, u_1 ) } \quad \text{and} \quad d_2 := \frac d{ ( d, u_2 ) }, \quad h_2 := \frac{u_2}{ ( d, u_2 ) }. \]
In order to evaluate these sums, we encode the additive twists by means of Dirichlet characters,
\[ e\left( \frac{y r_1 h_2}{d_2} \right) = \frac1{ \varphi(d_2) } \sum_{ \chi_2 \bmod d_2 } \overline{\chi_2}(y r_1 h_2) \tau(\chi_2), \]
so that we get
\[ M_{3 \text b}(d) = \frac{ \varphi(d) }{ \varphi(d_1) \varphi(d_2) } \sum_\thrln{ \chi_1 \bmod d_1 }{ \chi_2 \bmod d_2 }{ \chi_2 \equiv \overline{\chi_1} \bmod d } \overline{\chi_1}(-h_1) \overline{\chi_2}(h_2) \tau(\chi_1) \tau(\chi_2) W_{2, 1} W_{1, 2}, \]
where
\[ W_{i, j} := \sum_{ \left( r, u_i d \right) = 1 } f_{3, i}\left( \frac{h r}{u_i} \right) \psi_{\delta_i}( r u_i^\ast ) \overline{\chi_j}(r). \]

We use Mellin inversion to write these sums as follows,
\[ W_{i, j} = \frac1{2\pi i} \int_{ (\sigma) } \! \frac{ {u_i}^s }{ h^s } \hat f_{3, i}(s) Z_3(s) \, ds, \]
where \( \sigma > 1 \), and where
\[ \hat f_{3, i}(s) := \int \! f_{3, i}(\xi) \xi^{s - 1} \, d\xi \quad \text{and} \quad Z_3(s) := \sum_{ \left( r, \frac{h d}{u_i} \right) = 1 } \frac{ \psi_{\delta_i}( r u_i^\ast ) \overline{\chi_j}(r) }{ r^s }. \]
The Dirichlet series \( Z_3(s) \) converges absolutely for \( \Re(s) > 1 \), but it can be checked that an analytic continuation is given by
\[ Z_3(s) = \psi_{\delta_i}(u_i^\ast) \prod_{ p \mid u_i d } \left( 1 - \frac{ \overline{\chi_j}(p) }{ p^s } \right) \prod_{ p \mid u_i u_i^\ast d } \left( 1 - \frac{ \overline{\chi_j}(p) }{ p^{ 1 + s + \delta_i } } \right)^{-1} \frac{ L\left( s , \overline{\chi_j} \right) }{ L\left( 1 + s + \delta_i, \overline{\chi_j} \right) }. \]
We want to move the line of integration to \( \sigma = \varepsilon \), and the only pole we have to take care of lies at \( s = 1 \) and appears exactly when \( \chi_j \) is the principal character, in which case
\[ \Res{ s = 1 } \,\, Z_3(s) = \frac1{ \zeta(2 + \delta_i) } \frac{ \psi_{\delta_i}(u_i^\ast) \psi_0( u_i d ) }{ \psi_{1 + \delta_i}( u_i u_i^\ast d ) }. \]
Furthermore we have the following bound for \( \hat f_{3, i}(s) \),
\[ \hat f_{3, i}(s) \ll {A_1}^{ \Re(s) } \min\left\{ 1, \frac1{ |s| }, \frac1{ \Omega |s| |s + 1| } \right\}, \]
which also holds for its derivative with respect to \( \delta_i \), as well as the following bounds for the involved \(L\)-functions,
\[ \zeta(s) \ll |t|^{ \frac{1 - \sigma}2 + \varepsilon } \quad \text{and} \quad L\left( s, \overline{\chi_j} \right) \ll ( |t| d_j )^{ \frac{1 - \sigma}2 + \varepsilon } \]
(see \cite[(3)]{Kol79} for the latter).
When \( \chi_j \) is principal, we get this way an asymptotic formula for \( W_{i, j} \) of the form
\[ W_{i, j} = \frac{ \hat f_{3, i}(1) }{ \zeta(2 + \delta_i) } \frac{u_i}h \frac{ \psi_{\delta_i}(u_i^\ast) \psi_0( u_i d ) }{ \psi_{1 + \delta_i}( u_i u_i^\ast d ) } + \BigO{ \frac{x^\varepsilon}{\Omega^\frac12} }, \]
while otherwise we get an upper bound,
\[ W_{i, j} \ll x^\varepsilon \frac{ {d_j}^\frac12 }{\Omega^\frac12}. \]

Eventually this leads to
\begin{multline*}
  M_3 := \iiint \! \Delta_{\delta_1}(\xi_1 + h) \Delta_{\delta_2}(\xi_2 + h) F(\xi_1, \xi_2, \eta) C_3(h) \, d\eta d\xi_1 d\xi_2 \, d\eta \\
    + \BigO{ \frac{ x^{1 + \varepsilon} }{\Omega^\frac12} + \frac{ A_1 x^{1 + \varepsilon} }{D_0} + \frac{ x^{1 + \varepsilon} D_0}{A_1 \Omega} },
\end{multline*}
with \( F(\xi_1, \xi_2, \eta) \) as defined in \eqref{eqn: definition of F(xi_1, xi_2, eta)}, and with
\begin{multline*}
  C_3(h) := \frac1{ \zeta(2 + \delta_1) \zeta(2 + \delta_2) } \sum_\twoln{h_1, h_2 \mid h}{ u_1 \mid h_1, \,\, u_2 \mid h_2 } \frac{ {h_1}^{\delta_1} {h_2}^{\delta_2} u_1 u_2 }{ h^{ 2 + \delta_1 + \delta_2 } } \psi_{\delta_1}\left( \frac{h_1}{u_1} \right) \psi_{\delta_1}\left( \frac{h_2}{u_2} \right) \\
    \cdot \sum_{ \left( d, \frac{h_1 h_2}{u_1 u_2} \right) = 1 } \frac1{ d^{1 + \delta_1 + \delta_2} } \frac{ \mu\left( \frac d{ (d, u_1) } \right) \mu\left( \frac d{ (d, u_2) } \right) }{ \varphi\left( \frac d{ (d, u_1) } \right) \varphi\left( \frac d{ (d, u_2) } \right) } \frac{ \psi_0(d) \psi_0( u_1 d ) \psi_0( u_2 d ) }{ \psi_{1 + \delta_1}( h_1 d ) \psi_{1 + \delta_2}( h_2 d ) }.
\end{multline*}
One can easily check that \( C_3(1) = C_{\delta_1, \delta_2} \), with \( C_{\delta_1, \delta_2} \) as defined in \eqref{eqn: defintion of C(h)}, and that
\[ \gamma_3(h) := \frac{ C_3(h) }{ C_3(1) } \]
is a multiplicative function in \(h\).
A much more tedious calculation then shows that \( \gamma_3(h) \) and \( \gamma_{\delta_1, \delta_2}(h) \) indeed agree on prime powers and hence must be the same function.
Clearly, the optimal value for \( D_0 \) is
\[ D_0 = A_1 \Omega^\frac12, \]
and we finally get \eqref{eqn: evaluation of Sigma_3}.

\bibliography{The_shifted_convolution_of_generalized_divisor_functions}

\begin{thebibliography}{10}

\bibitem{Blo08}
V.~Blomer.
\newblock The average value of divisor sums in arithmetic progressions.
\newblock {\em Q. J. Math.}, 59(3):275--286, 2008.

\bibitem{BV87}
V.~A. Bykovski{\u\i} and A.~I. Vinogradov.
\newblock Inhomogeneous convolutions.
\newblock {\em Zap. Nauchn. Sem. Leningrad. Otdel. Mat. Inst. Steklov. (LOMI)},
  160(Anal. Teor. Chisel i Teor. Funktsii. 8):16--30, 296, 1987.

\bibitem{DI82}
J.-M. Deshouillers and H.~Iwaniec.
\newblock Kloosterman sums and {F}ourier coefficients of cusp forms.
\newblock {\em Invent. Math.}, 70(2):219--288, 1982/83.

\bibitem{Dra15}
S.~Drappeau.
\newblock Sums of {K}loosterman sums in arithmetic progressions, and the error
  term in the dispersion method.
\newblock {\em arXiv:1504.05549v3 [math.NT]}, 2015.
\newblock arXiv preprint.

\bibitem{FT85}
{\'E}.~Fouvry and G.~Tenenbaum.
\newblock Sur la corr\'elation des fonctions de {P}iltz.
\newblock {\em Rev. Mat. Iberoamericana}, 1(3):43--54, 1985.

\bibitem{Hea86}
D.~R. Heath-Brown.
\newblock The divisor function {$d_3(n)$} in arithmetic progressions.
\newblock {\em Acta Arith.}, 47(1):29--56, 1986.

\bibitem{Kim03}
H.~H. Kim.
\newblock Functoriality for the exterior square of {${\rm GL}_4$} and the
  symmetric fourth of {${\rm GL}_2$}.
\newblock {\em J. Amer. Math. Soc.}, 16(1):139--183, 2003.
\newblock With appendix 1 by Dinakar Ramakrishnan and appendix 2 by Kim and
  Peter Sarnak.

\bibitem{Kol79}
G.~Kolesnik.
\newblock On the order of {D}irichlet {$L$}-functions.
\newblock {\em Pacific J. Math.}, 82(2):479--484, 1979.

\bibitem{Lin63}
J.~V. Linnik.
\newblock {\em The dispersion method in binary additive problems}.
\newblock Translated by S. Schuur. American Mathematical Society, Providence,
  R.I., 1963.

\bibitem{Meu01}
T.~Meurman.
\newblock On the binary additive divisor problem.
\newblock In {\em Number theory ({T}urku, 1999)}, pages 223--246. de Gruyter,
  Berlin, 2001.

\bibitem{Mot80}
Y.~Motohashi.
\newblock An asymptotic series for an additive divisor problem.
\newblock {\em Math. Z.}, 170(1):43--63, 1980.

\bibitem{Mot94}
Y.~Motohashi.
\newblock The binary additive divisor problem.
\newblock {\em Ann. Sci. \'Ecole Norm. Sup. (4)}, 27(5):529--572, 1994.

\bibitem{Top15}
B.~Topacogullari.
\newblock On a certain additive divisor problem.
\newblock {\em arXiv:1512.05770v2 [math.NT]}, 2015.
\newblock arXiv preprint.

\bibitem{Top16}
B.~Topacogullari.
\newblock The shifted convolution of divisor functions.
\newblock {\em Q. J. Math.}, Advance Access published May 3, 2016.
\newblock doi:10.1093/qmath/haw010.

\end{thebibliography}

\end{document}